\documentclass{amsart}
\usepackage{txfonts}
\newcommand{\bt}{\begin{Theorem}}
\newcommand{\et}{\end{Theorem}}
\newcommand{\bi}{\begin{itemize}}
\newcommand{\ei}{\end{itemize}}
\newcommand{\bea}{\begin{eqnarray}}
\newcommand{\ba}{\begin{array}}
\newcommand{\eea}{\end{eqnarray}}
\newcommand{\ea}{\end{array}}

\def\g{\mathfrak{g}}
\def\v{\mathfrak{v}}
\def\z{\mathfrak{z}}

\newcommand{\what}{\widehat}

\newtheorem{Definition}{Definition}[section]
\newtheorem{Theorem}[Definition]{Theorem}
\newtheorem{Lemma}[Definition]{Lemma}

\newtheorem{Remark}[Definition]{Remark}

\newcommand{\be}{\begin{equation}}
\newcommand{\ee}{\end{equation}}
\newcommand{\bvb}{\begin{verbatim}}
\newcommand{\evb}{\end{verbatim}}
\newcommand{\Bea}{\begin{eqnarray*}}
\newcommand{\Eea}{\end{eqnarray*}}
\newcommand{\newsection}{\setcounter{equation}{0}}

\newcommand{\R}{\mathbb R}%
\newcommand{\C}{\mathbb C}%
\newcommand{\N}{\mathbb N}%

\begin{document}

\title [Beurling's theorem and
$L^p-L^q$ Morgan's theorem]{Beurling's theorem and $L^p-L^q$
 Morgan's theorem for step two nilpotent Lie groups}
\author[Parui]{Sanjay\ Parui}
\address[Sanjay\ Parui]{Stat--Math Unit\\
Indian Statistical Institute\\
203 B.T. Road,Kolkata- 700108,India,  email: {\tt
sanjay\_r@isical.ac.in}}

\author[Sarkar]{Rudra P.\ Sarkar}
\address[Rudra P.  \ Sarkar]{Stat--Math Unit\\
Indian Statistical Institute\\
203, B.T. Road, Kolkata- 700108, India, email: {\tt
rudra@isical.ac.in}}

\subjclass[2000]{22E30, 43A80}

\keywords{Uncertainty principle, Beurling's theorem, Morgan's
theorem, Nilpotent Lie groups}

\footnote{The first author is supported by a research fellowship
of National Board for Higher Mathematics, India.}

\begin{abstract} We prove Beurling's theorem and
 $L^p-L^q$ Morgan's theorem for step two nilpotent Lie groups.
 These two theorems together imply a group of uncertainty
 theorems.

\end{abstract}
\maketitle

\section{Introduction}
\setcounter{equation}{0} Roughly speaking the {\em Uncertainty
Principle} says that ``A nonzero function $f$ and its Fourier
transform $\widehat f$ cannot be sharply localized
simultaneously". There are several ways of measuring localization
of a function and depending on it one can formulate different
versions of qualitative uncertainty principle (QUP). The most
remarkable result in this genre in recent times is due to
H\"{o}rmander \cite{HO} where decay has been measured in terms of
a single integral estimate involving  $f$ and $\widehat{f}$.
\begin{Theorem}{\em (H\"{o}rmander 1991)} Let $f\in L^2(\mathbb{R})$ be such that
$$\int_{\mathbb {R}}\int_{\mathbb R}|f(x)||\widehat f(y)|e^{|x||y|}~dx~dy<\infty.$$
Then $f=0$ almost everywhere.
\end{Theorem}
H\"{o}rmander attributes this theorem to A. Beurling. The above
theorem of H\"{o}rmander was further generalized by Bonami et al
\cite{BDJ} which also accommodates the optimal point of this
trade-off between the function and its Fourier transform:
\begin{Theorem}
\label{beurcor}Let $f\in L^2(\R^n)$ be such that
$$\int_{\R^n}\int_{\R^n}
 \frac{|f(x)||\widehat
f(y)|e^{|x||y|}}{(1+|x|+|y|)^N}dx~dy <\infty$$ for some $N\geq 0$.
Then $f=0$ almost everywhere whenever $N\leq n$. If $N>n$, then
$f(x)=P(x)e^{-a|x|^2}$ where  $P$ is a polynomial with $\deg
P<\frac{(N-n)}{2}$ and $a>0$.
\end{Theorem}
Following H\"{o}rmander we will refer to the theorem above simply
as Beurling's theorem.

This theorem is described as {\em master theorem} by some authors
as theorems of Hardy, Cowling-Price and some versions of Morgan's
as well as  $L^p-L^q$ Morgan's follow from it. (See Theorem
\ref{compact-up} for precise statements of these theorems.)

There is some misunderstanding regarding the implication of
Beurling's theorem. However it was observed by Bonami et. al.
(\cite{BDJ}) that Beurling's theorem does not imply  Morgan's
theorem in its sharpest form. Indeed  Beurling's theorem (Theorem
\ref{beurcor}) together with $L^p-L^q$ Morgan's theorem (Theorem
\ref{compact-up} (v)) can claim to be the master theorem.  We can
summarize the relations between these theorems on $\R^n$ in the
following
diagram. 

$$\begin{array}{ccccccc}&\Rightarrow&\text{Hardy's}&\,\,\,\,|\,\,\,\,&\text{Morgan's}&\Leftarrow&\\
 \text{Beurling's}&&\Uparrow&|&\Uparrow&&L^p-L^q\text{Morgan's}\\
&\Rightarrow&\text{Cowling-Price}&\,\,\,\,|\,\,\,\,&\text{Gelfand-Shilov}&\Leftarrow&\end{array}$$

The aim of this paper is to prove analogues of Beurling's theorem
and $L^p-L^q$ Morgan's theorem (Theorem \ref{beurcor}, Theorem
\ref{compact-up} (case v)) for the step two nilpotent Lie groups.
It is clear from the diagram that all other theorems mentioned
above follow from these two theorems. Note that the diagram above
remains unchanged when $\R^n$ is substituted by the step two
nilpotent Lie groups.

For the convenience of the  presentation and easy readability we
will first deal with the special case of the Heisenberg groups and
then extend the argument for general step two nilpotent Lie
groups. The organization of the paper is as follows. In section 2
we prove modified versions of Theorem \ref{beurcor} and Theorem
\ref{compact-up} for $\R^n$ which are important steps towards
proving those theorems for the class of groups mentioned above. In
section 3 we establish the preliminaries of the Heisenberg group
and prove the  two theorems for this group. In section 4 we put
the required preliminaries for general step two nilpotent Lie
groups. Finally in section 5 we prove the analogues of Beurling's
and $L^p-L^q$-Morgan's theorems for step two nilpotent groups. We
indicate how the other theorems of this genre follow from those
two theorems. We also show the necessity and sharpness of the
estimates used in the two theorems.

Some of the other theorems, which follow from Beurling's and
$L^p-L^q$-Morgan's (Hardy's  and Cowling-Price to be more
specific) were proved independently on Heisenberg groups or
nilpotent Lie groups in recent years by many authors (see
\cite{ACBS, BR, BS, BSS, KK, R} etc.). However we may note that
these theorems were proved under some restrictions. But as
corollaries of the Beurling's and $L^p-L^q$-Morgan's theorem we
get exact analogues of these theorems. We include a precise
comparison with the earlier results in the last section.
%
For a general survey on uncertainty principles on
different groups we refer to \cite{FS, T3}.

Acknowledgement: We thank the referee for many suggestions and
criticisms which improved the exposition.

\section{Euclidean Spaces}
\setcounter{equation}{0}
 We can state a group of uncertainty
principles in a compact form as follows:
\begin{Theorem}\label{compact-up}
Let $f$ be a measurable function  on $\mathbb{R}$. Suppose  for
some $a,b>0$, $p,q\in [1,\infty]$, $\alpha\geq 2$ and $\beta>0$
with $1/\alpha+1/\beta=1$, $f$ satisfies
$$e^{a|x|^\alpha}f\in L^p(\mathbb{R})~~ \mbox{and}~~ e^{b|y|^\beta}\widehat f\in L^q(\mathbb{R}).$$
If moreover
\begin{equation}(a\alpha)^{1/\alpha}(b\beta)^{1/\beta}>
(\sin\frac{\pi}{2}(\beta-1))^{1/\beta} \label{morgan-condtion}
\end{equation}
 then $f=0$ almost everywhere.
\end{Theorem}
The case
\begin{enumerate}
\item[(i)]
 $\alpha=\beta=2$ and $p=q=\infty$ is  Hardy's theorem.\\
\item[(ii)]
 $\alpha=\beta=2$  is Cowling--Price theorem.\\
\item[(iii)]
$\alpha>2$, $p=q=\infty$ is Morgan's theorem.\\
\item[(iv)] $\alpha>2$ and $p=q=1$ is Gelfand-Shilov theorem.\\
\item[(v)] $\alpha>2$, ~$p,q\in [1,\infty]$ is $L^p-L^q$ Morgan's
theorem.
\end{enumerate}
This theorem has ready generalization for $\R^n$ where by $|x|$ we
mean the Euclidean norm of $x$.

It is clear that we have two separate sets of results in the
theorem above namely the cases (i) and (ii) where $\alpha=2$ and
cases (iii), (iv), (v) where $\alpha>2$. Note that for the first
set, condition (\ref{morgan-condtion}) reduces to $ab>1/4$. Back
in 1934 Morgan \cite{M} observed that at the {\em optimal point}
of (\ref{morgan-condtion}) these  two sets behave differently. To
emphasize this we  consider  cases (i) and (iii) of Theorem 2.1 as
representatives of the two sets of results. It is known that when
$(a\alpha)^{1/\alpha}(b\beta)^{1/\beta}=
(\sin\frac{\pi}{2}(\beta-1))^{1/\beta}$ then in case (i) above $f$
is a constant multiple of the Gaussian. In {\em great contrast}
(see \cite{M}) there are uncountably many functions which satisfy
the estimates in case (iii) when
$(a\alpha)^{1/\alpha}(b\beta)^{1/\beta}=
(\sin\frac{\pi}{2}(\beta-1))^{1/\beta}$.

\subsection{Modified version of the Beurling's theorem:}
We will  state and prove
 a modified version of  Theorem \ref{beurcor}.
We need the following preparations.
Let $S^{n-1}$ denote the unit sphere in $\R^n$. For a suitable
function $g$ on $\R^n$, the  Radon transform $Rg$ is a function on
$S^{n-1}\times \R$, defined by
\begin{equation}Rg(\omega, r)=R_\omega
g(r)=\int_{x \cdot\omega=r}g(x)\, d\sigma_x, \label{radon}
\end{equation} where
$d\sigma_x$ denotes the $(n-1)$-dimensional Lebesgue measure on
the hyperplane $x\cdot\omega=r$ and $x\cdot\omega$ is the
canonical inner product of $x$ and $\omega$, i.e.,
$x.\omega=\sum_{i=1}^nx_i\omega_i$. Note that when $g\in
L^1(\R^n)$, then for any fixed $\omega\in S^{n-1}$, $Rg(\omega,
r)$ exists for almost every $r\in \R$ and is an $L^1$-function on
$\R$. It is also well known that (See \cite{Fol}, p. 185.)
\begin{equation}\widehat{R_\omega
g}(\lambda)=\widehat{g}(\lambda\omega).\label{radone}\end{equation}
Here $\widehat{R_\omega g}(\lambda)=\int_\R R_\omega
g(r)e^{-i\lambda r}dr$ and
$\widehat{g}(\lambda\omega)=\int_{\R^n}g(x)e^{-ix\cdot\lambda\omega}dx$.

We also need the following lemma:

\begin{Lemma}\label{lemma-quadratic} Let
$f_1(x)=P_1(x)e^{-\alpha_1x^2}$ and
$f_2(x)=P_2(x)e^{-\alpha_2x^2}$ be two functions on $\R$ where
$P_1, P_2$ are polynomials and $\alpha_1, \alpha_2$ are positive
constants. Suppose that
$$I_1=\int_\R\int_\R\frac{|f_1(x)||\widehat{f_2}(y)|e^{|xy|}Q(y)}{(1+|x|+|y|)^N}dx~dy<\infty$$
and
$$I_2=\int_\R\int_\R\frac{|f_2(x)||\widehat{f_1}(y)|e^{|xy|}Q(y)}{(1+|x|+|y|)^N}dx~dy<\infty$$
where $N$ is a positive integer and $Q$ is a polynomial. Then
$\alpha_1=\alpha_2$.
\end{Lemma}
\begin{proof}
We note that $\widehat{f_1}(y)=Q_1(y)e^{-\frac 1{4\alpha_1} y^2}$
and $\widehat{f_2}(y)=Q_2(y)e^{-\frac 1{4\alpha_2} y^2}$ where
$Q_1$ and $Q_2$ are polynomials with $\deg Q_1=\deg P_1$ and $\deg
Q_2=\deg P_2$. Then
\begin{eqnarray}I_1&=&\int_\R\int_\R\frac{e^{-\alpha_1 x^2+|xy|-\frac
1{4\alpha_2} y^2}Q(y)P_1(x)Q_2(y)}{(1+|x|+|y|)^N}dx~dy\nonumber\\
&=&\int_\R\int_\R\frac{e^{-(\sqrt{\alpha_1}|x|-\frac
1{2\sqrt{\alpha_2}}
|y|)^2}e^{(1-\frac{\sqrt\alpha_1}{\sqrt\alpha_2})|x||y|}Q(y)P_1(x)Q_2(y)}{(1+|x|+|y|)^N}dx~dy.\nonumber
\end{eqnarray}
Similarly  we get
\begin{eqnarray}I_2
&=&\int_\R\int_\R\frac{e^{-(\sqrt{\alpha_2}|x|-\frac
1{2\sqrt{\alpha_1}}
|y|)^2}e^{(1-\frac{\sqrt\alpha_2}{\sqrt\alpha_1})|x||y|}Q(y)P_2(x)Q_1(y)}{(1+|x|+|y|)^N}dx~dy.\nonumber
\end{eqnarray}
We fix an $\epsilon>0$ and consider the set $A_\epsilon=\{(x,y)\in
\R^2\mid |\sqrt\alpha_1|x|-\frac
1{\sqrt\alpha_2}|y||\le\epsilon\}$, which is clearly of infinite
measure Lebesgue measure.

If we assume that $\alpha_1<\alpha_2$, then
$\frac{\sqrt\alpha_1}{\sqrt\alpha_2}<1$ and hence there exists a
$M>0$ such that the integrand in $I_1$ is greater than $M$ on the
strip $A_\epsilon$. Therefore  $I_1=\infty$. This contradicts the
hypothesis that $I_1<\infty$. Similarly if we assume that
$\alpha_2<\alpha_1$, then $I_2=\infty$. This completes the proof.
\end{proof}

With this preparation we will now prove the following modified
Beurling's theorem for $\mathbb{R}^n$.
\begin{Theorem}\label{modified-beurling-for}
Suppose $f\in L^2(\R^n)$. Let for some $\delta>0$
\begin{equation}\int_{\R^n}\int_{\R^n}\frac{|f(x)||\widehat{f}(y)|e^{|x||y|}|Q(y)|^\delta}
{(1+|x|+|y|)^d}~dx~dy<\infty, \label{modified-bon}
\end{equation} where $Q$ is a polynomial of degree $m$. Then
$f(x)=P(x)e^{-a|x|^2}$ for some $a>0$ and polynomial $P$ with
$\deg P<\frac{d-n-m\delta}{2}$.
\end{Theorem}
\begin{proof}
\noindent{\bf Step 0:} As $\what{f}$ is not identically zero and
as $Q$ is a polynomial, the product $|\what{f}(y)| |Q(y)|^\delta$
is different from zero on a set of positive measure. Therefore  we
can assume that for some $y_0\in \R^n$,
$$\int_{\R^n}\frac{|f(x)|e^{|x||y_0|}}{(1+|x|+|y_0|)^d}dx<\infty.$$
As $f\in L^2(\R^n)$, it is a locally integrable function on $\R^n$
and hence for any $0<r<|y_0|$,
$\int_{\R^n}|f(x)|e^{r|x|}dx<\infty.$ This shows in particular
that $f\in L^1(\R^n)$. Indeed for the exponential weight
$e^{|y_0||x|}$ it is easy to see that $\what{f}$ is holomorphic in
a tubular neighbourhood in $\C^n$ around $\R^n$.

In (\ref{modified-bon}) we use polar coordinates for $y$, to see
that there exists a subset $S$ of $S^{n-1}$ with full surface
measure such that for every $\omega_2\in S$,
\begin{equation}
\int_{\R^n}\int_{\R}\frac{|f(x)||\widehat{f}(s\omega_2)||s|^{n-1}|Q(s\omega_2)|e^{|x||s|}}
{(1+|x|+|s|)^d}~ds~dx<\infty.
\end{equation} In view of (\ref{radone}) this is the same as for
every $\omega_2\in S$,
\begin{equation}
\int_{\R^n}\int_{\R}\frac{|f(x)||\widehat{R_{\omega_2}
f}(s)||s|^{n-1}|Q(s\omega_2)|e^{|x||s|}}
{(1+|x|+|s|)^d}~ds~dx<\infty\label{starting-point}.
\end{equation}

\noindent{\bf Step 1:} In this step we will show  that for any
$\omega_1\in S^{n-1}$ and
 $\omega_2\in S$,
\begin{equation}
\int_\R\int_\R\frac{R_{\omega_1}|f|(r)|\widehat{R_{\omega_2}f}(s)||s|^{n-1}|
Q(s\omega_2)|e^{|r||s|}} {(1+|r|+|s|)^d}~ds~dr<\infty.
\label{target0}
\end{equation}

We will break the above integral  into the following three parts
and show that each part is finite. That is we will show:
\begin{enumerate}
\item[(i)]
$$\int_{\R}\int_{|s|>L}\frac{ R(|f|)(\omega_1,r)
|\widehat{R_{\omega_2} f}(s)|e^{|r||s|}|s|^{n-1}|Q(s\omega_2)|}
{(1+|r|+|s|)^d}~ds~dr<\infty$$ for $L>0$ such that $L^2+L>d$.

\item[(ii)]
$$\int_{|r|>M}\int_{|s|\leq L}
\frac{{ R}(|f|)(\omega_1,r)|\widehat{R_{\omega_2}
f}(s)|e^{|r||s|}|s|^{n-1}|Q(s\omega_2)|}
{(1+|r|+|s|)^d}~ds~dr<\infty$$ for $M=2(L+1)$ and $L$ as in (i).

\item[(iii)]
$$\int_{|r|\leq M}\int_{|s|\le L}\frac{{
R}(|f|)(\omega_1,r)|\widehat{R_{\omega_2} f}(s)|
e^{|r||s|}|s|^{n-1}|Q(s\omega_2)|}{(1+|r|+|s|)^d}~ds~dr<\infty$$
for $M,L$ used in (i) and (ii).
\end{enumerate}

\noindent{\em Proof of} (i): It is given that $L+L^2>d$. We will
show that for any $s$ such that $|s|\geq L$,
\begin{equation}\frac{e^{|s||x|}}{(1+|x|+|s|)^d}\geq
\frac{e^{|s||\langle x,\omega_1\rangle |}}{(1+|\langle
x,\omega_1\rangle |+|s|)^d}. \label{useful-inequality}
\end{equation}
Let $F(z)=\frac{e^{\alpha z}}{(1+\alpha+z)^d}$ for $\alpha>0$ and
$\alpha+\alpha^2>d$. Then $F'(z)>0$ for any $z\geq 0$. Therefore,
if $z_1\geq z_2\geq 0$, then
\begin{equation}\frac{e^{\alpha
z_1}}{(1+\alpha+z_1)^d}\geq \frac{e^{\alpha
z_2}}{(1+\alpha+z_2)^d}. \label{elementary}
\end{equation}
Note that $|x|\geq |\langle x, \omega_1\rangle|$ for all $x\in
\R^n$ and $\omega_1\in S^{n-1}$. Now take $z_1=|x|$ and
$z_2=|\langle x, \omega_1\rangle |$. Then $z_1\geq z_2\geq 0$. We
take $\alpha=|s|\geq L$ to get (\ref{useful-inequality}).

We start now from (\ref{starting-point}) and break it up as:
\begin{equation}\int_\R\int_{x\cdot\omega_1=r}\int_\R\frac{|f(x)||\widehat{R_{\omega_2} f}(s)
|e^{|x||s|}|s|^{n-1}|Q(s\omega_2)|} {(1+|x|+|s|)^d}ds~
d\sigma_1~dr<\infty,
\end{equation}where $d\sigma_1$ denotes the Lebesgue measure on the
hyper plane $\{x: x\cdot\omega_1=r\}$. We use the inequality
(\ref{useful-inequality}) to obtain:
\begin{equation}\int_\R\int_{x\cdot\omega_1=r}\int_{|s|>L}\frac{|f(x)||\widehat{R_{\omega_2}
 f}(s)
|e^{|\langle x, \omega_1\rangle ||s|}|s|^{n-1}|Q(s\omega_2)|}
{(1+|\langle x, \omega_1\rangle |+|s|)^d}~ds~d\sigma_1~dr<\infty.
\end{equation}
Now we put $\langle x, \omega_1\rangle =r$ in the above integral
and use the definition of Radon transform to obtain,

\begin{equation}\int_\R\int_{|s|>L}\frac{
R(|f|)(\omega_1,r)|\widehat{R_{\omega_2}
f}(s)|e^{|r||s|}|s|^{n-1}|Q(s\omega_2)|}
{(1+|r|+|s|)^d}~ds~dr<\infty. \label{required-later}
\end{equation}
This proves (i). \vspace{.08in}

\noindent{\em Proof of} (ii): Let
$$I_2=\int_{|r|>M}\int_{|s|\le L}\frac{
R(|f|)(\omega_1,r)|\widehat{R_{\omega_2}
f}(s)|e^{|r||s|}|s|^{n-1}|Q(s\omega_2)|}{(1+|r|+|s|)^d}ds~dr.$$ It
is clear that, \Bea I_2&\leq &C\int_{|r|>M} \frac{{
R}(|f|)(\omega_1,r)|e^{L|r|}}{(1+|r|)^d}dr\\&=&C\int_{|r|>M}
\int_{x\cdot\omega_1=r}
\frac{|f(x)|e^{L|r|}}{(1+|r|)^d}~d\sigma_1~dr\\&=&CI_3,
\mbox{say}.\Eea We will show that $I_3$ is finite for  $M=2(L+1)$.

We have already observed that $\widehat f$ is real analytic on
$\R^n$ and hence $\widehat f(y)\neq 0$ for almost every
$y\in\R^n$. Therefore, from (\ref{starting-point}) we can get a
$s_0\in \R$ with $|s_0|>2L$ such that:
$$\int_{\R^n}\frac{|f(x)|e^{|x||s_0|}}{(1+|x|+|s_0|)^d}dx<\infty.$$

That is,
$$\int_{\R}\int_{x\cdot\omega_1=r}\frac{|f(x)|e^{|x||s_0|}}
{(1+|x|+|s_0|)^d}~d\sigma_1 ~dr<\infty.$$ Notice that
$|s_0|+|s_0|^2>d$, since $|s_0|>2L$ and $L+L^2>d$. Now applying
the argument of case (i) (see (\ref{elementary})) to $|s_0|$ we
get:
$$\frac{e^{|x||s_0|}}{(1+|x|+|s_0|)^d}\geq
\frac{e^{|\langle x, \omega_1\rangle | |s_0|}}{(1+|\langle x,
\omega_1\rangle |+|s_0|)^d} $$ as $|\langle x,\omega_1\rangle
|\leq |x|$. Therefore, \Bea&&
\int_{|r|>M}\int_{x\cdot\omega_1=r}\frac{|f(x)|
e^{|r||s_0|}}{(1+|r|+|s_0|)^d}~d\sigma_1~dr \\&\leq&
\int_{|r|>M}\int_{x\cdot\omega_1=r}\frac{|f(x)| e^{|\langle x,
\omega_1\rangle ||s_0|}}{(1+|\langle x, \omega_1\rangle
|+|s_0|)^d}~d\sigma_1~dr<\infty\Eea from the above observation.
Note that $M+M^2>d$ as $M=2(L+1)$ and $L+L^2>d$. Applying the
argument of case (i) again (see (\ref{elementary})) this time with
$\alpha=|r|>M$ and $z_1=|s_0|, z_2=2L$ we get,
$$\frac{e^{|s_0||r|}}{(1+|s_0|+|r|)^d}\geq
\frac{e^{2L|r|}}{(1+2L+|r|)^d}.$$ Therefore,
$$\int_{|r|>M}\int_{x\cdot\omega_1=r}
\frac{|f(x)|e^{2L|r|}}{(1+|r|+2L)^d}~d\sigma_1~dr<\infty.$$ From
this it is easy to see that
$$\int_{|r|>M}\int_{x\cdot\omega_1=r}\frac{|f(x)|e^{L|r|}}{(1+|r|)^d}
~d\sigma_1~dr<\infty$$ and hence, $I_3<\infty$. This completes the
proof of (ii).

\noindent{\em Proof of} (iii): As the domain $[-M,M]\times[-L, L]$
is compact and as
$$\frac{|\widehat{R_{\omega_2} f}(s)|
e^{|r||s|}|s|^{n-1}|Q(s\omega_2)|}{(1+|r|+|s|)^d}$$ is continuous
in this domain, the integral is bounded by $C\int_{-M}^{M}R
|f|(\omega_1, r)dr$. Now recall that $f\in L^1(\R^n)$. Therefore,
\bea \int_{-M}^{M}R |f|(\omega_1, r)dr\nonumber&\leq&\int_\R R
|f|(\omega_1,r)dr\nonumber\\&=&\int_\R\int_{x\cdot\omega_1=r}|f(x)|d\sigma dr\nonumber\\
&=&\int_{\R^n}|f(x)|~dx<\infty. \eea Thus from (i), (ii) and (iii)
we obtain (\ref{target0}). This completes step 1.

\noindent{\bf Step 2:} Using $|{ R}_\omega f (r)|\leq { R}_\omega
|f|(r)$ we see from (\ref{target0})  that for almost every
$\omega\in S^{n-1}$,
\begin{equation}
\int_\R\int_\R\frac{|R_{\omega}f(r)||\widehat{R_{\omega}f}(s)||s|^{n-1}|Q(s\omega)|^\delta
e^{|r||s|}} {(1+|r|+|s|)^d}~dr~ds<\infty.
\end{equation}
Now as for fixed $\omega$, $|s|^{n-1}|Q(s\omega)|^\delta$ is a
proper map in $s$ and as $R_\omega f$ as well as $\what{R_\omega
f}$ are locally integrable functions  we can apply  the
$1$-dimensional case of Theorem 1.2 to conclude that $R_\omega
f(r)=A_\omega(r)e^{-\alpha r^2}$, for some polynomial $A_\omega$
which depends on $\omega$ with $\deg
A_\omega<\frac{d-m\delta-n}{2}$ and $\alpha$ is a positive
constant. A priori, $\alpha$ also should depend on $\omega$. But
we will see below that $\alpha$ is actually independent of
$\omega$. It is clear that $\widehat{R_\omega
f}(s)=P_\omega(s)e^{-\frac 1{4\alpha}s^2}$, where $\deg P_\omega$
is same as $\deg A_\omega $. Consider $\omega_1, \omega_2\in S$
with $\omega_1\neq \omega_2$ for which $R_{\omega_1}$,
$R_{\omega_2}$ satisfy (\ref{target0}).
From the
 argument above it follows that $R_{\omega_1}
f(r)=A_{\omega_1}(r)e^{-\alpha_1 r^2}$ and $\widehat{R_{\omega_2}
f}(s)=P_{\omega_2}(s)e^{-\frac 1{4\alpha_2}s^2}$ for some positive
constants $\alpha_1, \alpha_2$. Therefore by Lemma
\ref{lemma-quadratic}, $\alpha_1=\alpha_2=\alpha$, say and
$\widehat{R_\omega f}(s)=P_{\omega}(s)e^{-\frac {1}{4\alpha}s^2}$.
\vspace{.2in}

 \noindent{\bf Step 3:} We will show that
$P_\omega(s)=P( s\omega)$ is a polynomial in $s\omega$, that is
$P$ is a polynomial in $\R^n$. We recall that $\widehat{R_\omega
f}(s)=\widehat{f}(s\omega)$ is a holomorphic function in a
neighbourhood around $0$ (see Step 0). We can write
$P_\omega(s)=\widehat{f}(s\omega)e^{\frac 1{4\alpha}
s^2}=\widehat{f}(s\omega)e^{\frac 1{4\alpha}
|s\omega|^2}=F(s\omega),$ say.

We write $F(s\omega)=\sum_{j=0}^ka_j(\omega)s^j$, where
$k=\max\limits_{\omega\in S^{n-1}} \deg
P_\omega<\frac{d-m\delta-n}{2}$.  Then for $j=0,1,\dots, k$
$$\left. \frac 1{j!}\, \frac{d^j}{ds^j}F(s\omega)\right |_{s=0}=a_j(\omega).$$
The left hand side is the restriction of a homogenous polynomial
of degree $j$ to $S^{n-1}$. Therefore $F(s\omega)$ is a polynomial
of degree $\leq k$ in $\R^n$. Therefore
$\widehat{f}(x)=P(x)e^{-\frac 1{4\alpha}|x|^2}$, where $\deg
P<\frac{d-m\delta-n}{2}$.
\end{proof}

\subsection{Modified version of the $L^p-L^q$ Morgan's theorem }
We will state and prove a modified version of $L^p-L^q$ Morgan's
theorem on $\mathbb{R}^n$.
\begin{Theorem}\label{modified-morgan2}
Let $f$ be a measurable function on $\mathbb{R}^n$. Suppose for
some $a,b>0$, $p,q\in[1,\infty]$, $\alpha>2$ and $\beta$ with
$1/\alpha+1/\beta=1$, $f$ satisfies the following conditions:
\begin{enumerate}
\item[(i)]$\int_{\mathbb {R}^n}e^{pa|x|^\alpha}|f(x)|^p~dx<\infty$,\\
\item[(ii)] $\int_{\mathbb {R}^n}e^{qb|y|^\beta}|\widehat
f(y)|^q|Q(y)|^\delta~dy<\infty$, where $Q(y)$ is a polynomial in
$y$ of degree $k$ and $\delta>0$.
\end{enumerate}
If
$(a\alpha)^{1/\alpha}(b\beta)^{1/\beta}>(\sin\frac{\pi}{2}(\beta-1))^{1/\beta}$
then $f=0$ almost everywhere. If
$(a\alpha)^{1/\alpha}(b\beta)^{1/\beta}=(\sin\frac{\pi}{2}(\beta-1))^{1/\beta}$
then there are infinitely many linearly independent functions
which satisfy {\em (i)} and {\em (ii)}.
\end{Theorem}
\begin{proof}
First we will see that the theorem is true for $n=1$. From the
hypothesis (i) it is clear that $f \in L^1(\R)$ and hence
$\widehat f$ is continuous. Also as $|Q(y)|^\delta$ is a proper
map, we immediately get
$$\int_{\mathbb{R}}e^{qb|y|^\beta}|\widehat
f(y)|^q~dy<\infty.$$

That is, $f$ satisfies all the hypothesis of Theorem
\ref{compact-up} case (v) and hence the theorem  for $n=1$
follows.

Now we assume that $n\geq 2$. Let us consider the case $p=q=1$ for
the sake of simplicity. For each $\omega\in s^{n-1}$
 \bea
\int_{\mathbb{R}}e^{a|r|^\alpha}|R_\omega f(r)|~dr
&\leq& \int_{\mathbb{R}}e^{a|r|^\alpha}R_\omega|f|(r)~dr\\
&=&\int_{\mathbb{R}}\int_{x\cdot\omega=r}|f(x)|~d\sigma~dr\nonumber\\
&\leq & \int_{\mathbb{R}}\int_{x\cdot\omega=r}e^{a|x|^\alpha}|f(x)|~d\sigma~dr\nonumber\\
&=&\int_{\mathbb{R}^n}e^{a|x|^\alpha}|f(x)|~dx<\infty\nonumber.
\eea Here $d\sigma$ denotes the measure on the hyperplane
$\{x:x\cdot\omega=r\}$. Using the polar coordinates we get \Bea
&&\int_{\mathbb{R}}\int_{S^{n-1}}e^{b|r|^\beta}|\widehat{R_\omega
f}(r)||r|^{n-1}|Q(r\omega)|^\delta~d\omega~dr\\
&=&\int_{\mathbb{R}}\int_{S^{n-1}}e^{b|r|^\beta}|\widehat
f(r\omega)||r|^{n-1}|Q(r\omega)|^\delta~d\omega~dr\\
&=&2\int_{\mathbb{R}^n}e^{b|y|^\beta}
|f(y)||Q(y)|^\delta~dy<\infty. \Eea
 Hence almost every $\omega\in S^{n-1}$
\begin{equation}
\int_{\mathbb{R}}e^{b|r|^\beta}|\widehat{R_\omega
f}(r)||r|^{n-1}|Q(r\omega)|^\delta~dr=\int_{\mathbb{R}}e^{b|r|^\beta}|\widehat
f(r\omega)||r|^{n-1}|Q(r\omega)|^\delta~dr<\infty.
\end{equation}
We can now  apply the one-dimensional case of the theorem proved
above to the function $R_\omega f$ to conclude that for almost
every $\omega\in S^{n-1}$, $\widehat{R_\omega f}(r)=\widehat
f(r\omega)=0$ whenever
$(a\alpha)^{1/\alpha}(b\beta)^{1/\beta}>(\sin\frac{\pi}{2}(\beta-1))^{1/\beta}$
and hence $f=0$ almost everywhere.\\
Given $a,b>0$ with
$(a\alpha)^{1/\alpha}(b\beta)^{1/\beta}>(\sin\frac{\pi}{2}(\beta-1))^{1/\beta}$
we can always choose $a^\prime<a$, $b^\prime<b$ such that
$(a^\prime\alpha)^{1/\alpha}(b^\prime\beta)^{1/\beta}>
(\sin\frac{\pi}{2}(\beta-1))^{1/\beta}$. If $p,q>1$, using
H\"{o}lder inequality together with the given hypothesis we get
$$\int_{\mathbb{R}^n}e^{a^\prime|x|^\alpha}|f(x)|~dx<\infty~~\mbox{and}~~
\int_{\mathbb{R}^n}e^{b^\prime|y|^\beta}|\widehat
f(y)||Q(y)|^\delta~dy<\infty.$$ Hence the first part of the
theorem follows.

For the last part: Let us define the function $f$ by
$$f(x)=-i\int_{C}z^\nu e^{z^q-qAz|x|^2}~dz$$ where $q=\frac{\alpha}{\alpha-2}, A^\alpha=\frac{1}{4}((\alpha-2)a)^2,
\nu=\frac{2m+4-\alpha}{2(\alpha-2)}$, $m\in\mathbb{R}$ and $C$ is
a path lies in the half plane $\Im z>0$, and goes to infinity, in
the directions $\theta=\arg z=\pm \theta_\circ$, where
$\frac{\pi}{2}q<\theta_\circ<\frac{1}{2}\pi$. Then Ayadi et al.
\cite{AM}  shows
 with the help of Morgan's \cite{M} method that
 \bea
f=\mbox{O}(|x|^me^{-a|x|^\alpha})~\mbox{ and}~  \widehat
f=\mbox{O}(|y|^{m^\prime} e^{-b|y|^\beta})\label{example-Morgan},
~\mbox{where} ~m^\prime=\frac{2m+n(2-\alpha)}{(2\alpha-2)}.\eea We
will apply this result to construct functions satisfying the
equality cases of the hypothesis. Assume that the degree of the
polynomial $Q$ is $k$.

For $p=\infty, q=\infty$, we choose $m^\prime$ so that
 $m^\prime+k\delta<0$ and
 $m=\frac{2m^\prime+n(2-\beta)}{(2\beta-2)}<0$ equivalently we
 choose
  $m^\prime<\min\{-k\delta,-\frac{n}{2}(2-\beta)\}$ and
construct a function $f$ satisfying (\ref{example-Morgan}). This
$f$ will satisfy both the hypothesis. If $p\neq \infty, q=\infty$
we choose $m^\prime$ satisfying $m^\prime+k\delta<0$ and
$m=\frac{2m^\prime+n(2-\beta)}{(2\beta-2)}<-\frac{n}{p}$ which
holds if and only if
$m^\prime<\min\{-k\delta,-\frac{n}{p}(\beta-1)-\frac{n}{2}(2-\beta)\}$
and construct the required function. For $p\neq\infty$ and
$q\neq\infty$ we have to choose
$m^\prime<\min\{-\frac{n+k\delta}{q},
-\frac{n}{p}(\beta-1)-\frac{n}{2}(2-\beta)\}$.
\end{proof}

\section{Heisenberg groups}
\setcounter{equation}{0} Main results in this section are
analogues of Theorem \ref{beurcor} and Theorem 2.1 (v) for the
Heisenberg groups. Let us first recall some basic facts of the
Heisenberg groups. The $n$-dimensional Heisenberg group $H^n$ is
$\mathbb{C}^n\times \mathbb{R}$ equipped with the following group
law
 $$(z,t)(w,s)=(z+w,t+s+\frac{1}{2}\Im(z.\bar w)),$$ where $\Im(z)$ is the
 imaginary part of $z\in \C$.
 For each
$\lambda\in \R\setminus\{0\}$ there exists an irreducible unitary
representation $\pi_\lambda$ realized on $L^2(\R^n)$ given by
$$\pi_\lambda(z,t)\phi(\xi)=e^{i\lambda t}e^{i\lambda(x\cdot\xi+
\frac{1}{2}x\cdot y))}\phi(\xi+y),$$ for $\phi\in L^2(\R^n)$ and
$z=x+iy$. These are all the infinite dimensional irreducible
unitary representations of $H^n$ up to unitary equivalence. For
$f\in L^1(H^n)$, its group Fourier transform $\widehat f(\lambda)$
is defined by \bea\label{four}\widehat
f(\lambda)=\int_{H^n}f(z,t)\pi_\lambda(z,t)~dz~dt.\eea We define
$\pi_\lambda(z)=\pi_\lambda(z,0)$ so that
$\pi_\lambda(z,t)=e^{i\lambda t}\pi_\lambda(z,0)$. For $f\in
L^1(\C^n)$, we define the bounded operator $W_\lambda(f)$
 on $L^2(\R^n)$  by \bea
W_\lambda(f)\phi=\int_{\C^n}f(z)\pi_\lambda(z)\phi~dz~.\eea
 It is clear that
$\|W_\lambda(f)\|\leq \|f\|_1 $ and for $f\in L^1(\C^n)\cap
L^2(\C^n)$, it can be shown that $W_\lambda(f)$ is an
Hilbert-Schmidt operator and we have the Plancherel theorem
\bea\label{plancherel}
\|W_\lambda(f)\|_{\mbox{\tiny{HS}}}^2=(2\pi)^{n}|\lambda|^{-n}\int_{\C^n}|f(z)|^2~dz.
\eea Thus $W_\lambda$ is an isometric isomorphism between
$L^2(\C^n)$ and ${\mathcal{S}}_2$, the Hilbert space of all
Hilbert-Schmidt operators on $L^2(\R^n)$. This $W_\lambda(f)$ is
known as the Weyl transform of $f$. For $f\in L^1(H^n)$, let
$$f^\lambda(z)=\int_{-\infty}^{\infty}e^{i\lambda t}f(z,t)~dt$$  be the
inverse Fourier transform of $f$ in the $t$--variable. Then from
the definition of $\widehat f(\lambda)$, it follows that $\widehat
f(\lambda)=W_\lambda(f^\lambda)$. For $\lambda=1$ we define
$W(z)=W_1(z)$. For $x\in\R$ and $k\in\N$, the polynomial $H_k(x)$
of degree $k$ is defined by the formula \bea
H_k(x)=(-1)^ke^{x^2}\frac{d^k}{dx^k}(e^{-x^2}). \eea We define the
Hermite function $h_k(x)$ by \Bea h_k(x)=(2^k k!\sqrt
\pi)^{-1/2}H_k(x)e^{-\frac{x^2}{2}}. \Eea For
$\mu=(\mu_1,\cdots,\mu_n)\in \N^n$, the normalized Hermite
function $\Phi_{\mu}(x)$ on $\R^n$ is defined by \bea
\Phi_\mu(x)=h_{\mu_1}(x_1)\cdots h_{\mu_n}(x_n) .\eea  Hermite
functions are eigenfunctions of the Hermite operator
$H=-\triangle+|x|^2$ and they form an orthonormal basis for
$L^2(\R^n)$. Here $\triangle$ is the  Laplacian on $\mathbb{R}^n$.
For $\mu,\nu\in \N^n $, the special Hermite function
$\Phi_{\mu\nu}$ is defined by \bea
\Phi_{\mu,\nu}(z)=(2\pi)^{-\frac{n}{2}}\left(W(z)\Phi_{\mu},\Phi_\nu\right).
\eea These functions form an orthonormal basis for $L^2(\C^n)$ and
they are expressible in terms of Laguerre functions. For a
detailed account of Hermite and special Hermite functions we
refer to \cite{T5}.\\

With this preparation we will now  prove  a version of Theorem
\ref{beurcor} for $H^n$.
\begin{Theorem}\label{beurheisenb}
Suppose $f\in L^2(H^n)$ and for some $M,N\geq 0$, it satisfies
$$\int_{H^n}\int_{\R}\frac{|f(z,t)|\|\widehat f(\lambda)\|_{\mbox{\tiny{HS}}}e^{|t||\lambda|}}
{\left(1+|z|\right)^M\left(1+|t|+|\lambda|\right)^N}~|\lambda|^n~d\lambda~dz~dt<\infty.$$
Then
$f(z,t)=e^{-at^2}(1+|z|)^{M}\left(\sum\limits_{j=0}^m\psi_j(z)t^j\right)$,
where $\psi_j\in L^2(\C^n)$ and $m<\frac{N-n/2-1}{2}$.
\end{Theorem}
\begin{proof}
As in the case of $\R^n$, it can be verified that $f$ is
integrable in $t$-variable for almost every $z$.  For each pair
$(\phi,\psi)$, where $\phi,\psi\in L^2(\R^n)$ we consider the
function
$$F_{(\phi,\psi)}(t)=(2\pi)^{-\frac{n}{2}}\int_{\C^n}f(z,t)(1+|z|)^{-M}\overline{(W(z)\phi,\psi)}~dz.$$
Then it follows that
\bea\label{beurhei1}|\widehat{F_{(\phi,\psi)}}(\lambda)|&=&(2\pi)^{-\frac{n}{2}}
\left|\int_{\C^n}f^{-\lambda}(z)(1+|z|)^{-M}\overline{(W(z)\phi,\psi)}~dz\right|\\
&\leq&C\left(\int_{\C^n}|f^{-\lambda}(z)|^2~dz\right)^{1/2}\nonumber\\
&=& C |\lambda|^{n/2}\|\widehat
f(-\lambda)\|_{\mbox{\tiny{HS}}}\nonumber.\eea Therefore, from the
hypothesis we have
\Bea&&\int_{\R}\int_{\R}\frac{|F_{(\phi,\psi)}(t)||\widehat{F_{(\phi,\psi)}}(\lambda)|
e^{|t||\lambda|}|\lambda|^{n/2}}{\left(1+|t|+|\lambda|\right)^N}~dt~d\lambda\\
&\leq& C\int_{H^n}\int_{\R}\frac{|f(z,t)|\|\widehat
f(\lambda)\|_{\mbox{\tiny{HS}}}e^{|t||\lambda|}}
{\left(1+|z|\right)^{M}\left(1+|t|+|\lambda|\right)^N}|\lambda|^n~d\lambda~dz~dt<\infty.
\Eea
 Now applying Theorem \ref{modified-beurling-for} to the function
$F_{(\phi,\psi)}$ with $\delta=n/2$ we have
$F_{(\phi,\psi)}(t)=P_{(\phi,\psi)}(t)e^{-a(\phi,\psi)t^2}$, where
$P_{(\phi,\psi)}$ is a polynomial with $\deg<\frac{N-n/2-1}{2}$.
Keeping $\psi$ fixed, it can be shown that $a(\phi,\psi)=a(\psi)$
is independent of $\phi$. Similarly keeping $\phi$ fixed,  we can
show that $a(\phi,\psi)=a(\psi)=a$ is independent of
$(\phi,\psi)$. We recall that
$\{\Phi_{\alpha,\beta}:\alpha,\beta\in\N^n\}$ forms an orthonormal
basis for $L^2(\C^n)$. Now we take $\phi=\Phi_\alpha$ and
$\psi=\Phi_\beta$. Let
$F_{\alpha,\beta}=F_{(\Phi_\alpha,\Phi_\beta)}$ and
$P_{\alpha,\beta}=P_{(\Phi_\alpha,\Phi_\beta)}$. Since for each
$t\in\R$, $(1+|\cdot|)^{-M}f(\cdot,t)\in L^2(C^n)$, the sequence
$\{P_{\alpha,\beta}(t)\}\in l^2$ for all $t$. We write
$P_{\alpha,\beta}(t)=\sum\limits_{j=0}^m a_j(\alpha,\beta)t^j$,
$m<\frac{N-n/2-1}{2}$. Choose $t_i\in\R $ such that $t_i\neq t_j$,
for all $0\leq i,j\leq m$. We consider a system of linear
equations given by:
$$ \left(%
\begin{array}{cccc}
  1 & t_0 & \cdots & t_0^m \\
  1 & t_1 & \cdots & t_1^m \\
   \vdots&  \vdots &  \vdots &  \vdots \\
1 & t_m & \cdots & t_m^m \\
\end{array}%
\right)
\left(%
\begin{array}{c}
 \{ a_0(\alpha,\beta)\} \\
  \{a_1(\alpha,\beta)\}\\
  \vdots \\
 \{ a_m(\alpha,\beta)\} \\
\end{array}%
\right)
=\left(%
\begin{array}{c}
  \{P_{\alpha,\beta}(t_0)\} \\
  \{P_{\alpha,\beta}(t_1)\} \\
  \vdots \\
  \{P_{\alpha,\beta}(t_m)\} \\
\end{array}%
\right).$$ Since $t_i\neq t_j$ for all $ i\neq j$, the determinant
of the $(m+1)\times(m+1)$ Vandermonde matrix is nonzero.
Therefore, $\{a_j(\alpha,\beta)\}$ will be a linear combination of
members from $\{\{P_{\alpha,\beta}(t_j)\}: 0\leq j\leq m\}$
 and hence
$\{a_j(\alpha,\beta)\}\in l^2$ for each $0\leq j\leq m$. With this
observation we can write \Bea
(1+|z|)^{-M}f(z,t)&=&\left(\sum\limits_{\alpha,\beta}P_{\alpha,\beta}(t)\Phi_{\alpha,\beta}(z)\right)e^{-at^2}\\
&=&\left(\sum\limits_{\alpha,\beta}\left(\sum\limits_{j=0}^ma_j(\alpha,\beta)t^j\right)
\Phi_{\alpha,\beta}(z)\right)e^{-at^2}\\
&=&\left(\sum\limits_{j=0}^m\left(\sum\limits_{\alpha,\beta}a_j(\alpha,\beta)
\Phi_{\alpha,\beta}(z)\right)t^j\right)e^{-at^2}\\&=&\left(\sum\limits_{j=0}^m
\psi_j(z)t^j\right)e^{-at^2},\Eea where
$\psi_j(\cdot)=\sum\limits_{\alpha,\beta}a_j(\alpha,\beta)
\Phi_{\alpha,\beta}(\cdot) \in L^2(\C^n)$ .
\end{proof}

We will conclude this section by proving the following analogue of
$L^p-L^q$ Morgan's theorem for $H^n$.
\begin{Theorem}\label{Morgan-Heisen}
Suppose a function $f\in L^2(H^n)$  satisfies
\begin{enumerate}
\item[(i)] $\int_{H^n}e^{pa|(z,t)|^\alpha}|f(z,t)|^p~dz~dt<\infty$ and\\
\item[(ii)] $\int_{\mathbb{R}}e^{q|\lambda|^{\beta}}\|\widehat
f(\lambda)\|_{\tiny{\mbox{HS}}}^q|\lambda|^n~d\lambda<\infty$
\end{enumerate}
where $p,q\in[1,\infty]$, $a,b>0$, $\alpha>2$, $\beta>0$ and
$\frac{1}{\alpha}+\frac{1}{\beta}=1$.

If
$(a\alpha)^{1/\alpha}(b\beta)^{1/\beta}>(\sin\frac{\pi}{2}(\beta-1))^{1/\beta}$
then $f=0$ almost everywhere. But if
$(a\alpha)^{1/\alpha}(b\beta)^{1/\beta}=(\sin\frac{\pi}{2}(\beta-1))^{1/\beta}$
then there are infinitely many functions on $H^n$ satisfying {\em
(i)} and {\em (ii)}.

\end{Theorem}
\begin{proof}
First we note that $f\in L^1(H^n)$. We can choose $a^\prime<a$,
$b^\prime<b$ such that $(a^\prime\alpha)
^{1/\alpha}(b^\prime\beta)^{1/\beta}>(\sin\frac{\pi}{2}(\beta-1))^{1/\beta}$
and use H\"{o}lder's inequality  to show
\begin{enumerate}
\item[(i)$'$] $\int_{H^n}e^{a^\prime{|(z,t)|}^\alpha}|f(z,t)|~dz~dt<\infty$\\
\item[(ii)$'$]
$\int_{\mathbb{R}}e^{b^\prime|\lambda|^{\beta}}\|\widehat
f(\lambda)\| _{\tiny{\mbox{HS}}}|\lambda|^{n/2}~d\lambda<\infty$.
\end{enumerate} For each $(\mu,\nu)\in \mathbb {N}^n\times\mathbb{N}^n$, we define the
auxiliary function
$$F_{\mu,\nu}(t)=\int_{\mathbb{C}^n}f(z,t)\overline{\Phi_{\mu,\nu}(z)}~dz.$$
Using $(i)^\prime$ we have \Bea
&&\int_\mathbb{R}e^{a^\prime|t|^\alpha}|F_{\mu,\nu}(t)|~dt\\&=&\int_{\mathbb
{R}}
\int_{\mathbb{C}^n} e^{a^\prime |t|^\alpha}|f(z,t)||\Phi_{\mu,\nu}(z)|~dz\\
&\leq& \int_{\mathbb
{R}}\int_{\mathbb{C}^n}e^{a^\prime|(z,t)|^\alpha}|f(z,t)|~dz~dt\\&<&\infty.
\Eea On the other hand using $(ii)^\prime$ and the Plancherel
formula for the Weyl transform we have \Bea
&&\int_{\mathbb{R}}e^{b^\prime|\lambda|^\beta}|\widehat {F_{\mu,\nu}}(\lambda)|~d\lambda\\
&=&\int_{\mathbb{R}}e^{b^\prime|\lambda|^\beta}|\int_{\mathbb{C}^n}f^{-\lambda}(z)
\overline{\Phi_{\mu,\nu}(z)}~dz|d\lambda\\
&\leq& C \int_{\mathbb{R}}e^{b'|\lambda|^\beta}\|f^{-\lambda}(\cdot)\|_2~d\lambda\\
&=& \int_{\mathbb{R}}e^{b^\prime|\lambda|^\beta}\|\widehat
f(\lambda)\|_{\tiny{\mbox{HS}}}|\lambda|^{\frac{n}{2}}~d\lambda<\infty.
\Eea
 Applying  Theorem \ref{compact-up} (case (iv))  to the function $F_{\mu,\nu}$ we conclude
 that $F_{\mu,\nu}=0$ for every $(\mu,\nu)\in\mathbb{N}^n\times\mathbb{N}^n$.
 Since $\{\Phi_{\mu,\nu}:(\mu,\nu)\in\mathbb{N}^n\times\mathbb{N}^n\}$ form an
 orthonormal basis for $L^2(\mathbb{C}^n)$ we conclude that $f=0$
 almost everywhere.

Now for the second part of the theorem first we consider the case
$p=q=\infty$. We recall that  Morgan (see \cite{M}) constructed
enumerable examples of functions $h$ on $\mathbb{R}$ such that
$$h=O(|t|^me^{-a|t|^\alpha})~\mbox{ and}~ \widehat
h=O(|\lambda|^{m^\prime}e^{-b|\lambda|^\beta})$$ where
$m^\prime\in\mathbb{R}$ and
$m=\frac{2m^\prime-\beta+2}{2(\beta-1)}$.

We shall use these functions to construct required functions on
$H^n$. We define a function $f$ on $H^n$ as follows:
$f(z,t)=g(z)h(t)$, where $g$ is a smooth function on
$\mathbb{C}^n$ with compact support and $h$ is as above. Now it is
easy to see from the Plancherel formula for the Weyl transform
that
\begin{enumerate}
\item[(i)]$f=O(|(z,t)|^me^{-a|(z,t)|^\alpha})$\\
\item[(ii)] $O(|\lambda|^n\|\widehat
f(\cdot)\|_{\tiny\mbox{HS}})=O(|\lambda|^{n/2+m^\prime}e^{-b|\lambda|^\beta})$.
\end{enumerate}
 We choose $m^\prime<-\frac{n}{2}$ so that $m=\frac{2m^\prime-\beta+2}{2(\beta-1)}<0$. Then $f$ satisfies the required
estimates in the case $p=q=\infty$. For $p\neq\infty$ and
$q=\infty$ we choose $m^\prime<-\frac{n}{2}$ so that $m<
-\frac{(2n+1)}{p}$. If we choose
$m^\prime<\min\{-\frac{(n+1)}{q}+\frac{n}{2},-(2n+1)\frac{\beta-1}{p}+\frac{\beta-2}{2}\}$,
then $m<-\frac{(2n+1)}{p}$. With this choice of $m^\prime$ it is
easy to see that $f$ satisfies the required estimates with
$q\neq\infty$ and $p\in[1,\infty]$.
\end{proof}

\section{Step two Nilpotent Lie
Groups} \setcounter{equation}{0} Let $G$ be a step two connected
simply connected nilpotent Lie group. Then its Lie algebra
 $\g$ has the decomposition $\g=\v\oplus \z$, where $\z$ is the centre of
 $\mathfrak{g}$ and $\v$
is any subspace of $\g$ complementary to $\z$. We choose an inner
 product on $\g$ such that $\v$ and $\z$ are orthogonal. We fix an orthonormal basis
$\mathcal{B}=\{e_1,e_2\cdots,e_{m},T_1,\cdots,T_k\}$ so that
$\mathfrak {v}=\mbox{span}_{\mathbb{R}}\{e_1,e_2\cdots,e_m\}$ and
$\mathfrak{z}=\mbox{span}_{\mathbb{R}} \{T_1,\cdots,T_k\}$. Since
 $\g$ is nilpotent the exponential map is an analytic diffeomorphism. We can
 identify $G$ with $\v\oplus\z$ and write $(X+T)$ for
$\exp(X+T)$ and denote it by $(X,T)$
 where $X\in\v$ and $T\in\z$. The product law on $G$ is given by the
 Baker-Campbell-Hausdorff formula :
 $$(X,T)(X^\prime,T^\prime)=(X+X^\prime,T+T^\prime+\frac{1}{2}[X,X^\prime])$$ $\mbox{for all}~
 X,X^\prime\in\v$ and $T,T^\prime\in\z$.
\subsection{Representations of step two nilpotent Lie groups}\setcounter{equation}{0}
A complete account of representation theory for general connected
simply connected nilpotent Lie groups can be found in \cite{CG}.
Representations of step two connected simply connected nilpotent
groups the Plancherel theorem is described in \cite{R}.  We
briefly recall the basic facts to make this paper self contained.
Let $\mathfrak{g}^\ast$, $\mathfrak{z}^\ast$ be the real dual of
$\mathfrak{g}$ and $\mathfrak{z}$ respectively. For each
$\nu\in\mathfrak{z}^\ast$ consider the bilinear form $B_\nu$ on
$\mathfrak{\mathfrak{v}}$ defined by
$$B_\nu(X,Y)=\nu([X,Y])~\mbox{for all}~ X,Y\in\mathfrak{v}.$$ Let
$$\mathfrak{r}_\nu=\{X\in\mathfrak{v}
:\nu([X,Y])=0~\mbox{for all}~Y\in \mathfrak{v}\}.$$  Let $X_i=e_i$
for  $1\leq i\leq m$ and $X_{m+i}=T_i$ for  $1\leq i\leq k$. Then
$\mathcal{B}=\{X_1,\cdots,X_m,X_{m+1},\cdots,X_{m+k}\}$. Let
$\mathcal{B}^\ast=\{X_1^\ast,\cdots,X_m^\ast,X_{m+1}^\ast,\cdots,X_{m+k}^\ast\}$
be the dual basis of $\mathcal{B}$. Let $\mathfrak{m}_\nu$ be the
orthogonal complement of $\mathfrak{r}_\nu$ in $\mathfrak{v}$.
Then the set $ \mathcal
{U}=\{\nu\in\mathfrak{z}^\ast:\mbox{dim}~(\mathfrak{m}_\nu)~\mbox{is
maximum}\}$ is a Zariski open subset of $\mathfrak{z}^\ast$. Since
$B_\nu$ is an alternating bilinear form, $\nu\in \mathcal {U}$ has
an even number of jump indices independent of $\nu$. The set of
jump indices is denoted by $S=\{j_1,j_2\cdots,j_{2n}\}$.  Let
$T=\{n_1,n_2,\cdots,n_r,m+1,\cdots,m+k\} $ be the complement of
$S$ in $\{1,2,\cdots,m,m+1,\cdots,m+k\}$. Let
$$V_S={\mbox{span}}_{\R}\{X_{j_1},\cdots,X_{j_{2n}}\},$$
$$V_T={\mbox{span}}_{\R}\{X_{m+1},\cdots,X_{m+k},X_{n_i}:n_i\in
T\}~\mbox{and}~\widetilde {V}_T={\mbox{span}}_{\R}\{X_{n_i}:n_i\in
T\},$$
$$V_T^\ast={\mbox{span}}_{\R}\{X_{m+1}^\ast,\cdots,X_{m+k}^\ast,X_{n_i}^\ast:n_i\in
T\} ~\mbox{and} ~\widetilde
{V}_T^\ast={\mbox{span}}_{\R}\{X_{n_i}^\ast:n_i\in T\}.$$ The
irreducible unitary representations relevant to Plancherel measure
are parametrized by the set $\Lambda=\widetilde V_T^\ast\times
\mathcal{U}$.

If there exist $\nu\in\mathfrak{z}^\ast$ such that  $B_\nu$ is
nondegenerate then we call the group, a step two nilpotent group
with MW-- condition or step two MW group. In this case
$T=\{m+1,\cdots,m+k\}$ and $\mathcal{U}=\{\nu\in
\mathfrak{z}^\ast: B_\nu~\mbox{ is nondegeneate}\}$. The
irreducible unitary representations relevant to Plancherel measure
 will be parametrized by
$\Lambda=\{\nu\in\mathfrak{z}^\ast: B_\nu~\mbox{is
nondegenerate}\}$.

For
$$(X,T)=\mbox{exp}(\sum\limits_{j=1}^{m}x_jX_j+\sum\limits_{j=1}^kt_jX_{j+m}),~~x_j,t_j\in\R,$$
 we define its norm by $$|(X,T)|=(x_1^2 +\cdots + x_m^2+ t_1^2 +\cdots
+t_k^2)^{1/2}.$$ The map
\Bea(x_1,\cdots,x_m,t_1\cdots,t_k)\longrightarrow&\sum\limits_{j=1}^{m}
x_jX_j+\sum\limits_{j=1}^kt_jX_{j+m}\longrightarrow&\mbox{exp}\left(\sum\limits_{j=1}^{m}
x_jX_j+\sum\limits_{j=1}^kt_jX_{j+m}\right)\Eea takes Lebesgue
measure $dx_1\cdots dx_m dt_1\cdots dt_k$ of $\R^{m+k}$ to Haar
measure on $G$.  Any measurable function $f$ on $G$ will be
identified with a function on $\R^{m+k}$.  We identify
$\mathfrak{g}^\ast$ with $\R^{m+k}$ with respect to the basis
$\mathcal{B}^\ast$ and introduce the Euclidean norm relative to
this basis.
\subsubsection{Step two groups without MW--condition}
In this case $\mathfrak{r}_\nu\neq\{ 0$\} for each $\nu\in
\mathcal{U}$.  Then $B_\nu|_{\mathfrak{m}_\nu}$ is nondegenerate
and hence $\dim \mathfrak{m}_\nu$ is $2n$. From the properties of
an alternating bilinear form there exists an orthonormal basis
$$\{X_1(\nu),Y_1(\nu),\cdots,X_n(\nu),Y_n(\nu),Z_1(\nu),\cdots,
Z_r(\nu)\}$$ of $\mathfrak{v}$ and positive numbers $d_i(\nu)>0$
such that
\begin{enumerate}
\item [(i)]$\mathfrak{r}_\nu=\mbox{span}_{\R}~\{Z_1(\nu),\cdots,
Z_r(\nu)\} $,\\
\item[(ii)] $\nu([X_i(\nu),Y_j(\nu)])=\delta_{i,j}d_j(\nu),1\leq
i,j\leq n$.
\end{enumerate} We call
the basis
$$\{X_1(\nu),\cdots,X_n(\nu),Y_1(\nu),\cdots,Y_n(\nu),Z_1(\nu),\cdots,Z_r(\nu),T_1,
\cdots,T_k\}$$ almost symplectic basis. Let
$\xi_\nu={\mbox{span}}_{\R}\{X_1(\nu)\cdots,X_n(\nu)\}$ and
$\eta_\nu= \mbox{span}_{\R}~\{Y_1(\nu),\cdots,Y_n(\nu)\}$. Then we
have the decomposition
$\g=\xi_\nu\oplus\eta_\nu\oplus\mathfrak{r}_\nu\oplus\mathfrak z$.
We denote the element $\mbox{exp}(X+Y+Z+T)$ of $G$ by $(X,Y,Z,T)$
for $X\in\xi_\nu,Y\in\eta_\nu,Z\in \mathfrak{r}_\nu,T\in
\mathfrak{z}$. Further we can write
$$(X,Y,Z,T)=\sum\limits_{j=1}^n x_j(\nu)X_j(\nu)+\sum\limits_{j=1}^n y_j(\nu)Y_j(\nu)
+\sum\limits_{j=1}^r z_j(\nu)Z_j(\nu)+\sum\limits_{j=1}^k t_jT_j$$
and denote it by $(x,y,z,t)$ suppressing the dependence of $\nu$
which will be understood from the context. If we take
$\lambda\in\Lambda$ then it can be written as $\lambda=(\mu,\nu)$,
where $\mu\in \widetilde
V_T^\ast={\mbox{span}}_{\R}~\{X_{n_i}^\ast: 1\leq i\leq r\}$ and
$\nu\in\mathcal{U}$. Therefore,
$\lambda=(\mu,\nu)\equiv\sum\limits_{i=1}^r\mu_{i} X_{n_i}^\ast
+\sum\limits_{i=1}^m\nu_i T_i^\ast$. Let
$\lambda^\prime\in\mathfrak{g}^\ast$ such that
$\lambda^\prime(X_{j_i})=0$ for $1\leq i\leq 2n$ and the
restriction of $\lambda^\prime$ to $V_T^\ast$ is
$\lambda=(\mu,\nu)$. Let $\widetilde
\mu_i=\lambda^\prime(Z_i(\nu))$ and consider the map \bea
\label{phimap}A_\nu:\widetilde V_T^\ast\rightarrow
{\mbox{span}}_{\R}~\{Z_1(\nu)^\ast,\cdots,Z_r(\nu)^\ast\}\eea
given by $A_\nu(\mu_1,\cdots,\mu_r)=(\widetilde
\mu_1,\cdots,\widetilde\mu_r)$. Then it has been shown in \cite{R}
that $|\mbox{det}J_{A_\nu}|=\frac{\mbox{Pf}(\nu)}{d_1(\nu)\cdots
d_n(\nu)}$, where $J_{A_\nu}$ is the Jacobian matrix of $A_\nu$
and $\mbox{Pf}(\nu)$ is the Pfaffian of $\nu$. Consider the map
\bea D_\nu:\{X_{j1},\cdots, X_{j2n}\}\rightarrow
\{X_1(\nu),\cdots,X_n(\nu), Y_1(\nu),\cdots, Y_n(\nu)\}\eea then
it has been shown
$|\mbox{det}(J_{D_\nu})|=|\mbox{det}(J_{A_\nu})|^{-1}$ in
\cite{R}.

We take $\lambda=(\mu,\nu)\in\Lambda$.  Using the almost
symplectic basis we describe an irreducible unitary representation
$\pi_{\mu,\nu}$ of $G$ realized on $L^2(\eta_\nu)$ by the
following action:
\Bea\left(\pi_{\mu,\nu}(x,y,z,t)\phi\right)(\xi)&=&\mbox{exp}({i\sum\limits_{j=1}^k\nu_j
 t_j+i\sum\limits_{j=1}^r\widetilde{\mu}_jz_j
+i\sum\limits_{j=1}^nd_j(\nu)(x_j\xi_j+\frac{1}{2}x_jy_j)})\phi(\xi+y)\Eea
for all $\phi\in L^2(\eta_\nu)$.

We define the Fourier transform of $f\in L^1(G)$ by $$\widehat
f(\mu,\nu)=\int_{\mathfrak{z}}\int_{\mathfrak{r}_\nu}
\int_{\eta_\nu}\int_{\xi_\nu}
f(x,y,z,t)\pi_{\mu,\nu}(x,y,z,t)~dx~dy~dz~dt$$ for
$\lambda=(\mu,\nu)\in\Lambda$. For
$\widetilde\mu\in\mathfrak{r}_\nu^\ast,\nu\in\mathfrak{z}^\ast$ we
let
$$f^{\nu}(x,y,z)=\int_{\mathfrak{z}}\mbox{exp}({i\sum\limits_{j=1}^k\nu_jt_j})f(x,y,z,t)~dt~\mbox{and}$$
$$f^{\widetilde\mu,\nu}(x,y)=\int_{\mathfrak{r}_\nu}\int_{\mathfrak z}\mbox{exp}({i\sum
\limits_{j=1}^k \nu_j t_j+i\sum\limits_{j=1}^r\widetilde
\mu_jz_j})f(x,y,z,t)~dt~dz.$$  If $f\in L^1\cap L^2(G)$ then
$\widehat f(\mu,\nu)$ is an Hilbert--Schmidt operator and we have
(see \cite{R}) \bea (2\pi)^{-n}\prod_{j=1}^n d_j(\nu)\|\widehat
f({\mu},\nu)\|_{\mbox{\tiny{HS}}}^2=\int_{\eta_\nu}\int_{\xi_\nu}
|f^{\widetilde\mu,\nu}(x,y)|^2~dx~dy.\eea Now integrating both
sides on $\widetilde V_T^\ast$ with respect to the usual Lebesgue
measure on it and applying the transformation given by the
function $A_\nu$ in (\ref{phimap}) we get

\bea (2\pi)^{-(n+r)}\mbox{Pf}(\nu)\int_{\widetilde
V_T^\ast}\|\widehat f({\mu},\nu)\|_{\mbox{\tiny{HS}}}^2
 d{\mu}&=&(2\pi)^{-r}\int_{\mathfrak{r}_\nu^\ast}
\int_{\eta_\nu}\int_{\xi_\nu}|f^{\widetilde\mu,\nu}(x,y)|^2~dx~dy~d\widetilde\mu\nonumber\\
&=&\int_{\mathfrak{r}_\nu}
\int_{\eta_\nu}\int_{\xi_\nu}|f^{\nu}(x,y,z)|^2~dx~dy~dz\nonumber\\
&=&\int_{\mathfrak{v}}|f^\nu(x,y,z)|^2~dx~dy~dz.\nonumber\eea The
Plancherel formula takes the following form:
$$(2\pi)^{-(n+r+k)}\int_\Lambda\|\widehat f(\mu,\nu)\|_{\mbox{\tiny{HS}}}^2 ~\mbox{Pf}(\nu)~d\mu ~d\nu=
\int_{G}|f(x,y,z,t)|^2~dx~dy~dz~dt$$ which holds for all
$L^2$-functions by density argument.

\subsubsection{Step two MW groups }
In this case the representations are parametrized by the Zariski
open set $\Lambda=\{\nu\in\mathfrak{z}^\ast:B_\nu~\mbox{is
nondegenerate}\}$ and is given by:
 \bea
(\pi_\nu(x,y,t)\phi)(\xi)\nonumber&=&\mbox{exp}({i\sum\limits_{j=1}^k\nu_j
t_j+i \sum\limits_{j=1}^nd_j(\nu)
(x_j\xi_j+\frac{1}{2}x_jy_j)})\phi(\xi+y)\eea for all $\phi\in
L^2(\eta_\nu)$. In this case
$\mbox{Pf}(\nu)=\prod_{j=1}^nd_j(\nu)$. We define the Fourier
transform of $f\in L^1(G)$ by
$$\widehat f(\nu)=\int_{\mathfrak z}\int_{\eta_\nu}\int_{\xi_\nu}f(x,y,t)
\pi_\nu(x,y,t)~dx~dy~dt$$ for all $\nu\in\Lambda$. We also define
$$f^\nu(x,y)=\int_{\mathfrak z}\mbox{exp}({i\sum\limits_{j=1}^k\nu_j
t_j})f(x,y,t)~dx~dy~dt$$ for all $\nu\in \Lambda$.  If $f\in
L^1\cap L^2(G)$ then $\widehat f(\nu)$ is an Hilbert-Schmidt
operator and \bea\label{planchMW} \mbox{Pf}(\nu)\|\widehat
f(\nu)\|_{\mbox{\tiny{HS}}}^2&=&(2\pi)^n\int_{\eta_\nu}\int_{\xi_\nu}|f^\nu(x,y)|^2~dx~dy
=(2\pi)^n\int_{\mathfrak{v}}|f^\nu(x,y)|^2~dx~dy \nonumber.\eea
The Plancherel formula takes the following form: \bea
\label{planch-MW} (2\pi)^{-(n+k)}\int_{\Lambda}\|\widehat
f(\nu)\|_{\mbox{\tiny{HS}}}^2~\mbox{Pf}(\nu)~d\nu=\int_{G}|f(x,y,t)|^2~dx~dy~dt\eea
which  holds for all $L^2$-functions by density argument.

%
\section{Beurling's  and $L^p-L^q$-Morgan's theorem for step two Nilpotent Lie groups }\setcounter{equation}{0}
In what follows we will use the coordinates given by the following
basis of $\mathfrak g$.
$$\{X_{j_1},\cdots,X_{j_n}, X_{j_{n+1}},\cdots,X_{j_{2n}},X_{n_1},\cdots,X_{n_r},X_{m+1},\cdots,
X_{m+k}\}.$$   Precisely
$$(x,y,z,t)\equiv\sum\limits_{i=1}^nx_iX_{j_i}+\sum\limits_{i=1}^ny_iX_{j_{n+i}}+\sum\limits
_{i=1}^rz_iX_{n_i}+\sum\limits_{i=1}^kt_iX_{m+i}.$$

We shall first take up the following analogue of Beurling's
theorem for step two nilpotent groups.
\begin{Theorem}\label{for-step-2}Suppose $f\in L^2(G)$ and for some $M,N\geq
0$, it satisfies
\Bea&&\int_{\Lambda}\int_{\mathfrak{g}}\frac{|f(x,y,z,t)|\|\widehat
f(\mu,\nu)\|_{\tiny\mbox{HS}}
e^{|z||\mu|+|t||\nu|}}{(1+|(x,y)|)^M\left(1+|(z,t)|+|(\mu,\nu)|\right)^N}\\&&
\times|\mbox{Pf}(\nu)|~dx~dy~dz~dt~d\mu~d\nu<\infty.\Eea Then \Bea
&&f(x,y,z,t)\\&=&(1+|(x,y)|)^M\left(\sum\limits_{|\gamma|+|\delta|\leq
l}\Psi_{\gamma,\delta}(x,y)z^\gamma
t^\delta\right)e^{-a(|z|^2+|t|^2)} \Eea where
$\Psi_{\gamma,\delta}\in L^2(V_{S})$ and l is an nonnegative
integer.
\end{Theorem}
\begin{proof}
As in the case of $\R^n$ we can verify that $f$ is integrable in
$(z,t)$ for almost every $x, y$. For each Schwartz function $\Phi$
on $V_{S}$ let us consider the function $F_\Phi$ defined by \Bea
&&
F_\Phi(z,t)\\&=&\int_{V_{S}}f(x,y,z,t)(1+|(x,y)|)^{-M}\overline{\Phi(x,y)}~dx~dy.\Eea
It follows that
$$|F_\Phi(z,t)|\leq C \int_{V_S}|f(x,y,z,t)|~dx~dy.$$
 For all $(\mu,,\nu)\in\widetilde V_T^\ast\times\mathcal{U}$
$$\widehat{F_\Phi}(\mu,\nu)=\int_{V_S}f^{\mu,\nu}(x,y)(1+|(x,y)|)^{-M}
\overline{\Phi(x,y)}~dx~dy$$  where
$$f^{\mu,\nu}(x,y)=\int_{\widetilde V_T}e^{i\mu(z)+i\nu(t)}f(x,y,z,t)~dz~dt.$$ Using
Cauchy-Schwarz inequality, we get \Bea
&&|\widehat{F_\Phi}(\mu,\nu)|\\&\leq&C\left(\int_{V_S}|f^{\mu,\nu}(x,y)|^2~dx~dy\right)^{1/2}
\\&=&\left(\int_{V_S}\left|\int_{\widetilde V_T}e^{i\mu(z)}f^{\nu}(x,y,z)~dz\right|^2~dx~dy\right)^{1/2}.
 \Eea  Writing down the above integral with respect to almost
symplectic basis we have \Bea &&|\widehat{F_\Phi}(\mu,\nu)|\\
 &\leq&\left(\int_{\xi_\nu \oplus\eta_\nu}|f^{\widetilde{\mu},\nu}(x(\nu),y(\nu))|^2~dx(\nu)~dy(\nu)\right)^{1/2}\\
&=&(2\pi)^{-{n/2}}\left(\prod\limits_{j=1}^nd_j(\nu)\right)^{1/2}\|\widehat
f(\mu,\nu)\|_{\tiny\mbox{HS}}. \Eea Therefore, \Bea
&&\int_{\Lambda}\int_{\widetilde
V_T\oplus\mathfrak{z}}\frac{|F_\Phi(z,t)||\widehat
F_\Phi(\mu,\nu)|
e^{|\mu||z|+|t|\nu|}}{\left(1+|(z,t)|+|(\mu,\nu)|\right)^N(1+\prod\limits_{j=1}
^nd_j(\nu))}\\
&&\times|\mbox{Pf}(\nu)|~dz~dt~d\mu~d\nu\\
&\leq&C\int_{\Lambda}\int_{V_{S}\oplus\widetilde
V_T\oplus\mathfrak{z}}\frac{|f(x,y,z,t)|\|\widehat
f(\mu,\nu)\|_{\tiny{\mbox{HS}}}e^{|z||\mu|+|t||\nu|}}{(1+|(x,y)|)^M\left(1+|(z,t)|+|(\mu,\nu)|\right)^N}
\\&&\times|\mbox{Pf}(\nu)|~dx~dy~dz~dt~d\mu~d\nu\\
&=&\int_{\Lambda}\int_{\mathfrak{g}}\frac{|f(x,y,z,t)|\|\widehat
f(\mu,\nu)\|_{\mbox\tiny{{HS}}}e^{|z||\mu|+|t||\nu|}}{(1+|(x,y)|)^M\left(1+|(z,t)+|(\mu,\nu)|\right)^N}
\\&&\times|\mbox{Pf}(\nu)|~dx~dy~dz~dt~d\mu~d\nu\\
&<&\infty. \Eea Since $\mathcal{U}$ is a set of full measure on
$\mathfrak{z}^\ast, $ and $\mbox{Pf}(\nu)$,
$\prod\limits_{j=1}^nd_j(\nu)$ are polynomial in $\nu $ using
Theorem \ref{modified-beurling-for} we have for each Schwartz
function $\Phi$
$$F_{\Phi}(z,t)=P_{\Phi}(z,t)e^{-a(\Phi)|(z,t)|^2}$$
where $a(\Phi)>0$ and
$$P_{\Phi}(z,t)=\sum\limits_{|\gamma|+|\delta|\leq m}
a_{(\gamma,\delta)}(\Phi)z^\gamma t^\delta$$ and $m $ is
independent of $\Phi$. It is easy to see that $a(\Phi)=a$ is
independent of $\Phi$. Finally choosing $\Phi_\alpha$ from the
orthonormal basis $\{\Phi_\alpha(x,y):\alpha\in \mathbb{N}^{2n}\}$
for $L^2(V_S)$ we can show as in the proof of  Theorem
\ref{beurheisenb} that
$$f(x,y,z,t)=(1+|(x,y)|)^M\left(\sum\limits_{|\gamma|+|\delta|\leq m}
\Psi_{\gamma,\delta}(x,y) z^\gamma
t^\delta\right)e^{-a|(z,t)|^2},$$ where $\Psi_{\gamma,\delta}\in
L^2(V_S).$
\end{proof}

\noindent{\bf Consequences of Beurling's theorem:} Let us note the
following  consequences of Beurling's theorem.
\begin{Theorem}{\em (Morgan's theorem, weak version)} Let $f$ be a measurable function $G$. suppose for some
$a,b>0$, $\alpha\geq 2$, $\beta>0$
\begin{enumerate}
\item[(i)] $|f(x,y,z,t)|\leq C e^{-a|(x,y,z,t)|^\alpha}$\\
\item[(ii)] $|\mbox{Pf}(\nu)|^{1/2}\|\hat
f(\mu,\nu)\|_{\tiny\mbox{HS}}\leq C e^{-b|(\mu,\nu)|^\beta}$
\end{enumerate}
where $1/\alpha+1/\beta=1$ and
$(a\alpha)^{1/\alpha}(b\beta)^{1/\beta}\geq 1$. Then $f=0$ almost
everywhere unless $\alpha=\beta=2$ and $ab=1/4$ in which case
$f(x,y,z,t)=\psi(x,y)e^{-a|(z,t)|^2}$ for some $\psi\in L^2(V_S)$
and $|\psi(x,y)|\leq C e^{-a|(x,y)|^2}$.
\end{Theorem}

\begin{proof} It is clear from the hypothesis that $f\in L^2(G)$.
Since $\alpha\geq 2$ we have $|(x,y,z,t)|^\alpha\geq
|(x,y)|^\alpha+|(z,t)|^\alpha$. Therefore from hypothesis (i) we
get $|f(x,y,z,t)|\leq C e^{-a|(x,y)|^\alpha}e^{-a|(z,t)|^\alpha}.$
Now the theorem can be obtained from Theorem \ref{for-step-2} by
applying the inequality
$|\xi|^\alpha/\alpha+|\eta|^\beta/\beta\geq |\xi\eta|$ and using
the fact $e^{-a|(x,y)|^\alpha}\in L^1(V_S)$.
\end{proof}
In the proof above  we have used the fact that $\alpha\geq 2$ to
split the function $e^{-a|(x,y,z,t)|^\alpha}$ as a product of a
function in $ L^1(V_S)$ and $e^{-a|(z,t)|^\alpha}$. This motivates
us to formulate the following version of Morgan's theorem.
\begin{Theorem}
Let $f\in  L^2(G)$. Suppose for some $a,b>0$, $\alpha,\beta>0$
\begin{enumerate}
\item[(i)] $|f(x,y,z,t)|\leq g(x,y) e^{-a|(z,t)|^\alpha}$, $g\in L^1(V_S)$\\
\item[(ii)] $|\mbox{Pf}(\nu)|^{1/2}\|\hat
f(\mu,\nu)\|_{\tiny\mbox{HS}}\leq C e^{-b|(\mu,\nu)|^\beta}$
\end{enumerate}
where $1/\alpha+1/\beta=1$ and
$(a\alpha)^{1/\alpha}(b\beta)^{1/\beta}\geq 1$. Then $f=0$ almost
everywhere unless $\alpha=\beta=2$ and $ab=1/4$ in which case
$f(x,y,z,t)=\psi(x,y)e^{-a|(z,t)|^2}$ for some $\psi\in L^2(V_S)$
and $|\psi(x,y)|\leq g(x,y)$.
\end{Theorem}
\begin{Theorem}{\em(Cowling-Price)}
Suppose $f\in L^1\cap L^2(G)$ and it satisfies the following
conditions.
\begin{enumerate}
\item[(i)] $\int\limits_G
e^{pa|(x,y,z,t)|^2}|f(x,y,z,t)|^p~dx~dy~dz~dt<\infty$ and\\
\item[(ii)]$\int_{\Lambda}e^{bq|(\mu,\nu)|^2}\|\widehat
f(\mu,\nu)\|_{\mbox{\tiny{HS}}}^q|\mbox{Pf}(\nu)|~d\mu~d\nu<\infty.$
\end{enumerate}
Then for $ab\geq 1/4$ and $\min\{p,q\}<\infty$, $f=0$ almost
everywhere.
\end{Theorem}
\begin{proof}
Using H\"{o}lder's inequality we can find $M,N>0$ such that
\begin{enumerate}
\item[$(i)^\prime$]$\int_{V_S\oplus\widetilde
V_T\oplus\mathfrak{z}}\frac{e^{a|(z,t)|^2}|f(x,y,z,t)|}{(1+|(x,y)|)^M(1+|(z,t)|)^N}
~dx~dy~dz~dt<\infty$ and\\
\item[$(ii)^\prime$]$\int_{\Lambda}\frac{e^{b|(\mu,\nu)|^2}\|\widehat
{f}(\mu,\nu)\|_{\tiny{\mbox{HS}}}}{(1+|(\mu,\nu)|)^N}|\mbox{Pf}(\nu)|~d\mu~d\nu<\infty$.
\end{enumerate}
Therefore using Theorem \ref{for-step-2} we can conclude that
$f=0$ almost everywhere when $ab\geq 1/4$ and
$\min\{p,q\}<\infty$.
\end{proof}
\begin{Theorem}{\em(Hardy's theorem)}
Suppose $f$ is a measurable function on $G$ which satisfies the
following conditions:
\begin{enumerate}
\item[(i)]$|f(x,y,z,t)|\leq
g(x,y)(1+|(z,t)|)^me^{-a|(z,t)|^2},\mbox{where}~g\in
L^1\cap L^2(V_S)$~and\\
\item[(ii)]$|\mbox{Pf}(\nu)|^{1/2}\|\widehat
f(\mu,\nu)\|_{\tiny{\mbox{HS}}}\leq
(1+|(\mu,\nu)|)^me^{-b|(\mu,\nu)|^2}.$
\end{enumerate}
Then $f=0$ almost everywhere if $ab>1/4$ and if $ab=1/4$ then
$f(x,y,z,t)=P(x,y,z,t)e^{-a|(z,t)|^2}$, where
$P(x,y,z,t)=\left(\sum\limits_{|\alpha|+|\delta|\leq
m}\psi_{\alpha,\delta}(x,y)z^\delta t^\alpha\right)$ and
$\psi_{\alpha,\delta}\in L^2(V_S)$
\end{Theorem}
We omit the proof which is a straight forward application of the
theorem above. \vspace{.2in}

\noindent{\bf Sharpness of the estimate in Beurling's theorem:} We
will show that the  condition used in Beurling's theorem is
optimal. For the sake of simplicity we consider the Heisenberg
group $H^n$. We suppose a function $f\in L^1\bigcap L^2(H^n)$
satisfies \bea \int_{H^n}\int_{\R}\frac{|f(z,t)|\|\widehat
f(\lambda)\|_{\mbox{\tiny{HS}}}e^{c|t||\lambda|}}
{\left(1+|z|\right)^M\left(1+|t|+|\lambda|\right)^N}~|\lambda|^n~d\lambda~dz~dt<\infty\label{1}
\eea for some $c>0$.
\begin{enumerate}
\item[(i)] If $c>1$, then $f$ satisfies the hypothesis of Theorem
3.1 and hence $f(z,t)=g(z)e^{-at^2}$ for some $g\in L^1\bigcap
L^2(\mathbb{C}^n)$ and $a>0$. Since by the Plancherel theorem
(\ref{plancherel}) $\|\widehat
f(\lambda)\|_{\tiny{\mbox{HS}}}=(2\pi)^{n/2}|\lambda|^{-n/2}\|g\|_2e^{-\frac{1}{4a}\lambda^2},$
it is easy to see that
$f$ cannot satisfy (\ref{1}) unless $f=0$ almost everywhere. \\
\item[(ii)] Now we suppose $c<1$. We choose $a,b>0$ such that
$ab=c^2$ and we construct the function $f(z,t)=g(z)P(t)e^{-a
t^2}$, where $g\in L^1\bigcap L^2(\mathbb{C}^n)$ and $P$ is a
polynomial of any degree. Then $f$ will satisfy  (\ref{1}).
Clearly for fixed $z\in \mathbb{C}^n$ these functions are linearly
independent in the variable $t$.
\end{enumerate}

We shall now prove an exact analogue of $L^p-L^q$-Morgan's theorem
for step two nilpotent Lie groups.
\begin{Theorem}Let $f\in  L^2(G)$.
Suppose for some $a,b>0$, $\alpha>2, \beta>0$
\begin{enumerate}
\item[(i)]$\int_{G}e^{pa|(x,y,z,t)|^\alpha}|f(x,y,z,t)|^p~dv~dt<\infty $ and \\
\item[(ii)]$\int_{\Lambda}e^{qb|(\mu,\nu)|^\beta}\|\widehat
f(\mu,\nu)\|_{\mbox{\tiny{HS}}}
^q|\mbox{Pf}(\nu)|~d\mu~d\nu<\infty,$
\end{enumerate}
where $1/\alpha+1/\beta=1$ and $p,q\in[1,\infty]$. Then $f=0$
almost everywhere whenever
$(a\alpha)^{1/\alpha}(b\beta)^{1/\beta}>(\sin\frac{\pi}{2}(\beta-1))^{1/\beta}$.
\end{Theorem}
\begin{proof}
As in the case of $\R^n$ we see that $f \in L^1(G)$.
 We note that it is sufficient to consider the case $p=q=1$ as in
the case of Heisenberg groups. Since $\prod\limits_{j=1}^n
d_j(\nu)$, $\mbox{Pf}(\nu)$ are polynomials in $\nu$, for any
$b^\prime<b$, applying Minkowski's integral inequality with
respect to the measures $dx~dy$ and
$d\mu~e^{b^\prime|\nu|^\beta}|\mbox{Pf}(\nu)| d\nu$ we get
\Bea&&\left(\int_{V_S}\left(\int_{\Lambda}
|f^{\mu,\nu}(x,y)|e^{b^\prime|\nu|^\beta}~d\mu~|\mbox{Pf}(\nu)|d\nu\right)^2
~dx~dy\right)^{1/2}\\
&\leq
&\int_{\Lambda}e^{b^\prime|\nu|^\beta}\left(\int_{V_S}|f^{\mu,\nu}(x,y)|^2~dx~dy\right)
^{1/2}~d\mu~|\mbox{Pf}(\nu)|d\nu\\
&\leq&
\int_{\Lambda}e^{b^\prime|(\mu,\nu)|^\beta}\left(\int_{V_S}|f^{\mu,
\nu}(x,y)|^2~dx~dy\right)^{1/2}~d\mu~|\mbox{Pf}(\nu)|d\nu\\
&=&C\int_{\Lambda}e^{b^\prime|(\mu,\nu)|^\beta}\|\widehat
f(\mu,\nu)\|_{\mbox{\tiny{HS}}}\left(\prod\limits_{j=1}^nd_j(\nu)\right)^{1/2}~d\mu~
|\mbox{Pf}(\nu)|d\nu\\
&\leq& \int_{\Lambda}e^{b|(\mu,\nu)|^\beta}\|\widehat
f(\mu,\nu)\|_{\mbox{\tiny{HS}}}|\mbox{Pf}(\nu)|d\mu~d\nu
\\&<&\infty.\Eea

This implies that for almost every $(x,y)\in V_S$
\bea\int_{\Lambda}e^{b^\prime|\nu|^\beta}|f^{\mu,\nu}(x,y)||\mbox{Pf}(\nu)|~d\mu~d\nu<\infty.
\eea Since $\Lambda=\mathcal U\times \widetilde{V}_T^*$, it
follows that \bea\int_{\mathcal
U}e^{b^\prime|\nu|^\beta}|f^{\mu,\nu}(x,y)||\mbox{Pf}(\nu)|~d\nu<\infty
\eea for almost every $(x, y)\in V_S$ and $\mu\in
\widetilde{V}_T^*$.
 From the hypothesis (i) with $p=1$, it is easy to see that for
 almost every $(x,y)\in V_S$
 $$\int_{\mathfrak{z}}\int_{\widetilde V_T}e^{a|(z,t)|^\alpha}|f(x,y,z,t)|~dz~dt<\infty.$$
Therefore for almost every $(x,y)\in V_S$
\bea&&\int_{\mathfrak{z}}e^{a|t|^\alpha}|f^{\mu}(x,y,t)|~dt\\
&\leq&\int_{\mathfrak{z}}\int_{\widetilde
{V_T}}e^{a|(z,t)|^\alpha}
f(x,y,z,t)~dz~dt\nonumber\\&<&\infty\nonumber\eea where
$f^\mu(x,y,t)=\int_{\widetilde {V_T}}e^{\mu(z)}f(x,y,z,t)~dz$.

As $\mathcal U$ is a set of full measure, we can now apply Theorem
\ref{modified-morgan2} to the function $f^\mu(x,y,t)$ to conclude
that for almost every $(x,y)\in V_S$, $f(x,y,z,t)=0$ whenever
$(a\alpha)^{1/\alpha}(b^\prime\beta)^{1/\beta}>(\sin\frac{\pi}{2}(\beta-1))^{1/\beta}$.
Since given $a,b>0$ with
$(a\alpha)^{1/\alpha}(b\beta)^{1/\beta}>(\sin\frac{\pi}{2}(\beta-1))^{1/\beta}$,
it is always possible to choose $b^\prime<b$ satisfying
$(a\alpha)^{1/\alpha}(b^\prime\beta)^{1/\beta}>(\sin\frac{\pi}{2}(\beta-1))^{1/\beta}$,
the theorem follows.
\end{proof}

\begin{Remark}
\end{Remark}
\begin{enumerate}
\item The Gelfand-Shilov theorem and the Morgan's theorem (in
their sharpest forms) are particular cases of Theorem 5.6 ($p=q=1$
and $p=q=\infty$ respectively) and thus are accommodated in that
theorem.

\item In Theorem 3.2 we have seen  example of functions which
satisfy the hypothesis  with
$(a\alpha)^{1/\alpha}(b\beta)^{1/\beta}=(\sin\frac{\pi}{2}(\beta-1))^{1/\beta}$
in the case of $H^n$. Similar construction can be carried out in
this case also.

\item All the theorems proved above for step two groups without MW
condition can be formulated and proved for step two MW groups with
obvious and routine modifications.

\end{enumerate}

\section{Concluding Remarks}\newsection
There are a few attempts (see \cite{KK,BS,BSS}) in recent times to
prove theorems of this genre  for general nilpotent Lie groups.
The basic step in these works is to build a new function on the
central variable satisfying the equivalent conditions. But in the
process the sharpness of the result is lost and hence it is not
possible to get the case of optimality.  It is unlikely that the
method pursued in those papers will generalize to the case of all
nilpotent Lie groups since the explicit formula for $\|\hat
f(\lambda)\|_{\mbox{\tiny HS}}$ is crucial
 in the proof of Beurling's theorem, which is unavailable in this generality.
  We refer to
the remark in \cite[p.~493]{KK} in this context. Our aim  in this
paper is to get the analogues of these two theorems discussed
above  which can accommodate the case of optimality and without
any restriction on $(p,q)$ and $(\alpha, \beta)$. This is the
reason we restrict ourselves to the  step two nilpotent Lie
groups.

We conclude the paper with a brief discussion on comparison with
the existing results of this genre.
Beurling's theorem, i.e. analogues of  Theorem 1.2 is not
considered  so far for any nilpotent Lie groups.  In \cite{BSS} an
analogue of Theorem 1.1 is proved for the special class of
nilpotent Lie groups of the form $\R^n \rJoin \R$. In \cite{PT} a
version of Theorem 1.1 is formulated for stratified step two
nilpotent Lie groups where the estimate involves the matrix
coefficients of the Fourier transform, instead of the operator
valued Fourier transform, which seems to be more restrictive. None
of these  theorems accommodated the optimal case of Beurling's
theorem. In contrast our theorem takes care of the optimal case
and is also a suitable version to get back the other QUP-results
in full generality in the context of step two nilpotent groups.

A version of the $L^p-L^q$ Morgan's theorem is proved in \cite{F}
only for Heisenberg groups.

As mentioned in the introduction, other theorems of this genre
which follow from either Beurling's or $L^p-L^q$-Morgan's theorem
were proved independently by many authors in nilpotent Lie groups.
Nevertheless none of these works dealt with the characterization
of the optimal case. There are also some other restrictions on the
hypothesis. For instance in \cite{R} Ray proved  the
Cowling--Price theorem for step two nilpotent Lie groups without
MW-conditions with the assumption $1\leq p\leq\infty,q\geq 2$ and
$ab>1/4$. This result is generalized in \cite{BS}  for any
nilpotent Lie group with the restriction $2\leq p,q\leq\infty$ and
$ab>1/4$. As shown above we can have the Cowling-Price theorem
with the original condition $1\leq p,q\leq\infty$ and $ab\geq 1/4$
as a consequence of Beurling's theorem.

In \cite{R} Ray also proved a version of Morgan's theorem which is
similar to  Theorem 5.3 (in fact slightly weaker than Theorem 5.3)
and again can be obtained as a consequence of Beurling's theorem.
We recall that only a weak version of Morgan's theorem follows
from Beurling's theorem, while the actual Morgan's theorem follows
from $L^p-L^q$-Morgan's theorem (see Remark 5.7). In \cite{ACBS}
Astengo et. al. proved a version of Hardy's theorem where they put
condition on the operator norm of the Fourier transform, instead
of the usual pointwise estimate. We note that only by a slight
modification of our proof, a Beurling's theorem can be obtained
where Hilbert-Schmidt norm of the Fourier transform  is replaced
by its operator norm. (We formulated the theorem using
Hilbert-Schmidt because it appears to be more natural.) As a
consequence we can get the theorem in \cite{ACBS}.

Recently an analogue of Beurling's theorem is proved for
Riemannian symmetric spaces in \cite{SS}. Due to the structural
difference, the statement as well as the method of proving the
theorem is different and it involves decomposing the statement in
$K$-types and treating each component separately.

We summarize our aim in this paper as to obtain the most natural
analogue of Beurling's and $L^p-L^q$ Morgan's theorem for step two
nilpotent groups which can  accommodate the optimal case and from
which we can get back the strongest version of the other theorems
in this genre as consequences.

\end{document}